\documentclass[12pt,twoside,bezier,reqno]{amsart}
\pdfoutput=1
\usepackage{amsmath,amssymb,amsfonts,epsfig,graphicx,euscript}%
\usepackage[pdftex,bookmarks,bookmarksnumbered,linktocpage,pdfstartview=FitH]{hyperref}
\usepackage{amsmath}
\usepackage{amsfonts}
\usepackage{graphicx}
\usepackage{caption}
\usepackage{subcaption}
\usepackage{float}
\usepackage[all]{hypcap}
\numberwithin{equation}{section}
\usepackage{cite}
\usepackage{epsfig}
\usepackage{amsmath}
\usepackage{amsfonts}
\usepackage{calc}
\usepackage{amsmath,amsthm,amsfonts,amssymb,amscd,euscript}
\pdfoutput=1
\hypersetup{colorlinks,%
citecolor=red,%
filecolor=blue,%
linkcolor=blue,%
urlcolor=blue,%
pdftex}

\newlength{\spacer}
\newsavebox{\mybox}

\newcommand{\bse}{\begin{subequations}}
\newcommand{\ese}{\end{subequations}}
\newcommand{\be}{\begin{equation}}
\newcommand{\ee}{\end{equation}}
\newcommand{\bea}{\begin{eqnarray}}
\newcommand{\eea}{\end{eqnarray}}
\newcommand{\ba}{\begin{array}}
\newcommand{\ea}{\end{array}}

\input{epsf}
\newcommand{\filebegin}{\begin{document}}
\newcommand{\fileend}{\end{document}}
\makeatother
\def\thefootnote{}
\newcommand{\lo}{\longrightarrow}
\newcommand{\NMM}{\hspace*{2mm}}

\renewcommand{\baselinestretch}{1.1}

\input{epsf}
\renewcommand{\baselinestretch}{1.1}
\def\n{\noindent}%

\numberwithin{equation}{section}
\def\mapdown#1{\Big\downarrow\rlap
{$\vcenter{\hbox{$\scriptstyle#1$}}$}}
\newtheorem{theorem}{Theorem}[section]
\newtheorem{lemma}[theorem]{Lemma}
\newtheorem{proposition}[theorem]{Proposition}
\newtheorem{corollary}[theorem]{Corollary}

\theoremstyle{definition}
\newtheorem{definition}[theorem]{Definition}
\newtheorem{example}[theorem]{\sc Example}
\newtheorem{xca}[theorem]{Exercise}
\theoremstyle{remark}
\newtheorem{remark}[theorem]{Remark}

\begin{document}

\setcounter{page}{1} \noindent

\vspace*{2cm}
\begin{center}
{\bf\large Semi-Helices of Euclidean Spaces}
 \\[0.5cm]

{\bf A. Heydari$^{1,a}$ \footnote{1) a.heidar140@gmail.com}, S. Amiri-Sharifi$^{2,b}$\footnote{2) s$_{-}$amirisharifi@sbu.ac.ir}}\\%
\vspace*{0.4cm}
{\it {a. Department of Mathematics, Firuzabad Institute of Higher  Education, Firuzabad, Iran}}  \\%
{\it {b. Department of Physics, Shahid Beheshti University, G.C., Evin, Tehran 19839, Iran}}  \\%
\vspace*{1cm}
\end{center}%
\vspace*{0.5cm}
\begin{quotation}
\noindent
{\footnotesize
{\sc Abstract.}
We introduce semi-helix hyper surfaces of Euclidean spaces. We also provide a local characterization of how these semi-helices are constructed.}
\end{quotation}


\markboth 
{A. Heydari, S. Amiri-Sharifi}
 {Semi-Helices of Euclidean Spaces}


\section{Introduction}

\vskip 0.4 true cm

A helix manifold is defined by the property that the angle between the tangent space at each point with a fix direction, in the ambient euclidean space, is constant. The importance of helices comes from their vast application in nature, sciences and also engineering of mechanical tools.They arise in the field of computer aided design, computer graphics, simulation of kinematic motion or design of highways, the shape of DNA and carbon nanotubes \cite{AGK}. Helical structures often appear in fractal geometry \cite{YTFX}. 

 The application of constant angle hyper-surfaces in physics of liquid crystal has been studied in \cite{PA} by Cermelli and Scala. In \cite{GH} it is observed that shadow boundaries are related to helix submanifolds. Also in \cite{FJJL} and \cite{FM}, helix surfaces have been studied in non flat ambient spaces. The authors in \cite{JG} have proceeded with "helix  submanifolds" i.e. submanifolds with constant angle between their tangent spaces and a fixed direction $d$.
 
 In this paper we generalize the concept of $'Helix'$, in the sense that the angle between tangent spaces and a fixed direction $d$ can vary in a given sufficiently small interval $(\theta_0 - \varepsilon , \theta_0 + \varepsilon)$. 

 The paper is organized as follows: Section \ref{Def} contains definition of semi-helix hypersurfaces and in section \ref{Con} we present a full description of how this semi-helices are constructed.  Finally, in the last part we show that all helices are locally constructed in the same way. 
  
   Throughout this paper, we follow Scala and Ruiz \cite{JG}. As is customary, we consider $<.,.>$ as an inner product on $\mathbb{R}^{n}$ and 
   $T_{p} M$ is the tangent space at a point $p$ on the hypersurface $M$ in $\mathbb{R}^n$.

\vskip 0.6 true cm

\section{\bf {\bf \em{\bf Semi-helix hypersurfaces}}} \label{Def}

\vskip 0.4 true cm
 Let $d \in \mathbb{R}^n$ be any direction (i.e. a unitary vector) and let 
 $V \subset \mathbb{R}^n$ be a linear subspace. The angle 
 $\theta$ between $d$ and $V$ is the angle between the vectors $d$ and $\pi_{V}(d)$, where 
 $\pi_{V}: \mathbb{R}^n \longrightarrow V $  is defined to be the orthogonal projection onto $V$. That is to say $cos(\theta):= <d , \pi_{V}(d)>$.
 
 The definition of $semi-helices$ is given here:
\begin{definition}\label{semi-helix}
Let $M \subset \mathbb{R}^n$ be a hypersurface and $d \neq 0$ be a unitary fixed vector in $\mathbb{R}^n$. We say that $M$ is a $semi-helix \ \ w.r. \ to \ d $ and an initially chosen angle $\theta_0,$ if there exists a fixed real numaber $0 \leq \varepsilon < \frac{\pi}{2},$ such that for all $p \in M$, the angle between
 $d$ and $T_pM$, which we denote it by $\theta := \theta(p)$, belongs to the interval $(\theta_0 - \varepsilon , \theta_0 + \varepsilon)$.
\end{definition}
The semi-helix submanifold $M$ with domain of angle
 $(\theta_0- \varepsilon,\theta_0 + \varepsilon)$ will be denoted by 
 $M_{\varepsilon}(\theta_0)$. \\
 
 Note that if $\varepsilon = 0$ then $M_0(\theta_0)= M(\theta_0)$ is a helix (manifold) of angle $\theta_0$.
\begin{remark} 
Note that the angle $\theta$ comes from the angle between two "directions": $d$ and the normal to the tangent space at each point. Hence if the spin of $\theta$ is in the direct sense, $\theta$ is considered to be $positive$ and otherwise $\theta$ is $negative$.
\end{remark}


\begin{remark}
Suppose $M, \ \theta_0$ and $d$ are the same as in definition
 \ref{semi-helix} and also assume that $\xi:M \longrightarrow \nu(M)$ is a unit normal vector field, where $\nu(M)$ denotes the normal vector boundle. Then we can use
  $<d,\xi_p>=sin(\theta_p)$ rather than $cos(\theta_p)$ to identify the domain of angle $\theta_p$.
\end{remark}


Now we give a method to construct semi-helices in $\mathbb{R}^n$. Without loss of generality we can assume that $\theta_0$ is $zero$.

\vskip 0.6 true cm

\section{\bf {\bf \em{\bf Construction}}}\label{Con}

\begin{figure}[tb]
\includegraphics[width=100mm]{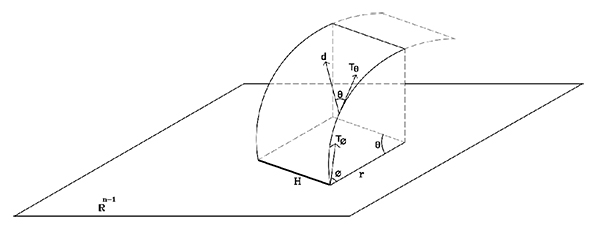}
\caption{}
\label{helix}
\end{figure}

\vskip 0.4 true cm
First we choose fixed 
$0 < \varepsilon < \frac{\pi}{2} $, and consider $(-\varepsilon, \varepsilon)$ the range of $\theta$.
 Let $H \subset \mathbb{R}^{n-1} $ be a hypersurface in $\mathbb{R}^{n-1}$
 and $\eta$ be a unitary normal vector field on $H$. We can immerse $H$ in $\mathbb{R}^n$ in a canonical way as $(x_1,x_2,...,x_{n-1})\longrightarrow (x_1,x_2,...,x_{n-1},0)$.

  without losing generality, we assume that $d=(0,...,0,1) \in \mathbb{R}^{n}$. put $M:= H \times S^{1}(r)$, in which 
  $S^{1}(r)$ is the circle of radius $r$. Now we define vector fields $T_{\theta}$ and $T_{\phi}$ as follows:
 
  \begin{equation}
   T_{\theta}(x) :=\sin(\theta)\eta(x) + \cos(\theta) d,
  \end{equation}
 \begin{equation}
   T_{\phi}(x) :=\cos(\phi)\eta(x) + \sin(\phi) d, 
  \end{equation}
where $x \in H$ and $\phi= \frac{\pi - \theta}{2}$. See figure \ref{helix}

Now we define the immersion 
$\mathit{f}_{\varepsilon}: H \times S^{1}(r) \longrightarrow \mathbb{R}^{n}$ as follows:

$$(x, (r \cos(\phi), r \sin(\phi))) \mapsto x + r \sqrt{(2(1- \cos \theta))} T_{\phi}(x) $$

Note that $- \varepsilon < \theta < \varepsilon$. Therefore for enough small $\varepsilon$, $\mathit{f}_{\varepsilon}$ is an immersion.
 
 
\begin{theorem}
The immersed submanifold
 $M_{\varepsilon} = \mathit{f}_{\varepsilon}(M)$ is a $semi-helix$.
\end{theorem}

\begin{proof}

First we will find a unitary normal vector field $\xi_{\theta}$ of the $M_{\varepsilon}$ and show that its angle with $d$ range in the interval $(-\varepsilon, \varepsilon)$. In order to do this we note that $f_{\varepsilon}|_{H} = Id_{H}$ and we can identify $H$ with $f_{\varepsilon}\{H \times {(r,0)}\}$ and $M$ with $M_{\varepsilon}$.Thus the tangent space of $M_{\varepsilon}$ at each point $f_{\varepsilon}(x,Y)$ is given by:
$$T_{(x,Y)}M_{\varepsilon}:= T_{x}H \oplus \mathbb{R} T_{\theta}(x), \  for \  any \  x\in H \  and \  Y \in S^1(r). $$
Now let $\xi_{\theta}(x):= - \cos(\theta) \eta(x)+ \sin(\theta) d$. Since $\eta$ and $d$ are orthogonal to $H$ and $<\xi_{\theta}(x),T_{\theta}(x)>=0$ thus $\xi_{\theta}(x)$ is a unitary normal vector field on $M_{\varepsilon}$. Now orthogonality of $\eta$ and $d$ imply that $<\xi_{\theta},d>=\sin(\theta)$ and this complete the proof.

\end{proof}



\vskip 0.6 true cm

\section{\bf {\bf \em{\bf Reconstruction}}}
\vskip 0.4 true cm

In this section we show that each $semi-helix \ \ M_{\varepsilon} \subset \mathbb{R}^{n}$, can be obtained locally in the form that described in the previous section. Suppose that $\xi_{\theta}$ and $T_{\theta}$ are the unitary normal and unitary tangent vector fields respectively on $semi-helix \ \ M_{\varepsilon}$ such that: 

$$d = \cos(\theta) T_{\theta}(p) + \sin(\theta)\xi_{\theta}(p)$$


\begin{lemma}\label{integral curve}
The integral curves of $T_{\theta}$ are sections of an $S^{1}(r)$ bundle over $M$.
\end{lemma}

\begin{proof}
Let $p \in M_{\varepsilon}$. by definition, in this point $\cos(\theta) \neq 0$. Now let $Q_{p}$ be a hyperplane that contains $p$ and is orthogonal to $d$. Put
 $H:= Q_{p} \cap M_{\varepsilon}$
 which is a submanifold of $M_{\varepsilon}$. There is a vector field $\eta$ which is normal to $H$ and satisfies:
 $T_{\theta}(p) = \sin(\theta) \eta(p) + \cos(\theta) d$.

   Now assume that $\alpha : (t_1 , t_2) \longrightarrow M_{\varepsilon}$
   is an integral curve of vector field $T_{\theta}$. So for $t_p \in (t_1, t_2)$ we have $\alpha(t_p)= p$ and
   
   $$\alpha^{'}(t) = T_{\theta}(\alpha(t)) = \sin(\theta) \eta(\alpha(t)) + \cos(\theta) d, \ \ \ \forall t \in (t_1, t_2)$$.

Since $\eta $ and $d $ are orthogonal, for all $t$ for which the angle $\theta$ corresponding with $\alpha(t) \in M_{\varepsilon}$ is positive, put:
$$e_{1}:= - \eta(\alpha(t)) \ \ , \ e_2 := d$$
 and for other $t$, put
 $$e_1:= \eta(\alpha(t)) \ \ ,\ e_2 := d$$
 
  Now we consider the extended basis $\beta= \{ e_1, e_2,...,e_n \}$ of $T_{\alpha(t)} \mathbb{R}^n$. For
   $\alpha^{'}(t) \in T_{\alpha(t)} \mathbb{R}^n$ we have: 
\begin{equation}\label{equation1}
    \alpha^{'}(t) = - \sin(\theta) e_{1} + \cos(\theta) e_{2}
\end{equation}

On the other hand consider the trivial chart 
$(Id , \mathbb{R}^n)$ around $\alpha(t)$, such that $\beta$ is a basis for $T_{\alpha(t)} \mathbb{R}^n \ w.r. \ $ to this chart. Let $\alpha(t)= (\alpha_{1}(t),..., \alpha_n(t))$ be the local representation of $\alpha \ \ w.r. \ $ to chart $(Id, \mathbb{R}^n)$. Thus we have: 
\begin{equation}\label{equation2}
\alpha^{'}(t) = \alpha_1^{'}(t) e_1 + ...+ \alpha_n^{'}(t) e_n.
\end{equation} 

equations \ref{equation1} and \ref{equation2} imply: 
$$\alpha_1^{'}(t) = - \sin(\theta) , \ \ \alpha_2^{'}(t) = \cos(\theta), \ \ \alpha_i^{'}(t)=0, \ \ (i=3,...,n)$$

Note that it is easy to show that, there is a linear relation between $\theta$ and $t$ of the form $\theta(t)=at + c$, where 
$a,c \in \mathbb{R}$ are constant. Therefor we have:
$$\alpha_1^{'}(t)= - \sin(at+c),$$ 
$$\alpha_2^{'}(t)= \cos(at+c),$$ 
$$\alpha_{i}^{'}(t) = 0, \ \ (i=3,...,n).$$ 

Thus:\\

 $$ \alpha_1(t)= \frac{1}{a} \cos(at+c) + c_1,$$
  $$\alpha_2(t)= \frac{1}{a} \sin(at+c) + c_2 ,$$
 $$ \alpha_{i}(t)=c_{i} \ \ \ (i=3,...,n).$$

where $c_{i} \ 's$ are constant. By using initial condition $\alpha(t_p)=p$, we can obtain the values of $c_{i} \ 's$. The last relation shows that $\alpha$ (the integral curve of $T_{\theta}$ ) is a section of a trivial $S^{1}(r)$ bundle on $M$ with 
$r= \frac{1}{a}$.
\end{proof}

\begin{theorem}

Each semi-helix hypersurface is locally isomorphic to the semi-helix hypersurface whose construction was explained in section \ref{Con}.
\end{theorem}

\begin{proof}
Let $M_{\varepsilon}$ be a semi-helix hypersurface $w.r. \ $ to a fixed direction $d$, and $T_{\theta}$ be the tangent component of  $d$ along $M_{\varepsilon}$ as mentioned in lemma 
\ref{integral curve}.

We choose an arbitrary $p \in M_{\varepsilon}$ and let $Q_{p}$ be a hyperplane passing through $p$ and orthogonal to $d$. Consider the submanifold $H:= Q_{p} \cap M_{\varepsilon} $ of $M_{\varepsilon}$. It is obvious that from each point $x \in H$, passes a curve $\alpha$ with values in $S^{1}(r)$  such that 
$T_{\theta}(x) = \alpha^{'}(t_x)$, $(\alpha(t_x)=x)$.

We denote the integral curve of $T_{\theta}$ passing through $p$ by $\alpha_{p}:I \longrightarrow M_{\varepsilon}$. In fact for all $t \in I$ the angle between $\alpha^{'}(t)$ and $d$ is $\theta_t$, that varies in the interval
 $(\varepsilon_1, \varepsilon_2)$. 
 
 Now we extend the curve $\alpha_p$ to
   $\beta_p : J \longrightarrow \mathbb{R}^n$ such that $\beta_p$ is a section of $S^{1}(r)$ bundle, 
 $I \subset J$ and $\beta_p | _{I} = \alpha_{p}$.
 Also there exist a point $q = \beta_{p}(t_{q}) \in \beta_{p}(J)$ such that the angle between $d$ and $\beta_{p}^{'}(t_q)$ is $zero$.
  
 Finally, we consider the hyperplane $Q_q$, that contains $q$ and is orthogonal to $d$. In fact $Q_q$ is a translation of 
  $Q_p$ and $q$ is a translation of $p$ under the same translation.
  
  We denote the translated submanifold $H$ of $Q_p$, by $H^{'}$. Let $\theta_p$ be the angle between $d$ and $T_{\theta}(p)$, then: 
  
$$p= q + r \sqrt{2(1- \cos(\theta_{p}))} T_{\phi_{p}}(q) $$ 
where $\phi_p = \frac{\pi - \theta_p}{2}$

 Now if $U \subset M_{\varepsilon}$ is a neighborhood around $p$, for all $y \in U \cap H$ we have: 
 $$y = y^{'} + r \sqrt{2(1- \cos(\theta_y))} T_{\phi_{y}}(y^{'})$$
 where $\phi_y = \frac{\pi - \theta_y}{2}$ and $y^{'} \in H^{'}$ is a translation of $y$. This indicates that 
 $U \subset \mathit{f}_{\varepsilon} (H^{'} \times S^{1}(r))$ and $f_{\varepsilon}(H^{'} \times S^{1}(r))$ is, as we showed in the previous section,a semi-helix hypersurface. This completes the proof.  
  
\end{proof}

\section{\bf {\bf \em{\bf Conclusion}}}

Our description indicates that semi-helix hypermanifolds of a euclidean space are locally isomorphic to the ones constructed in "Section 3", even though their global structure might involve some subtleties that could make them not to be isomorphic to our manifolds.

\section{\bf {\bf \em{\bf Acknowledgment}}}
The authors would like to thank Dr. M. Nadjafikhah and S. Shahriari for their kind and gracious help.

\vskip 2.1 true cm

\providecommand{\bysame}{\leavevmode\hbox
to3em{\hrulefill}\thinspace}



\begin{thebibliography}{99}

\bibitem{AGK} A. Jain,G. Wang,K. M. Vasquez, :{\em DNA triple helices: biological consequences therapeutic potential},Biochimie 90, 2008, pp. 1117-1130.

\bibitem{JG} A. J. Di Scala, G. Ruiz Hernandez,:{\em Helix submanifolds of euclidean spaces}, Monatsh. Math., vol. 157, 2009, pp.205-215.

\bibitem{BO} B. O'Neill,:{\em Semi-Reimannian Geometry with applications to relativity}, Academic press., New York, 1983.

\bibitem{FJJL} F. Dillen, J. Fastenakels, J. Van Der Veken, L. Vrancken,:{\em Constant Angle Surfaces in $\mathbb{S}^2 \times \mathbb{R}$}, Monatsh. Math., vol. 152, 2007, pp. 89-96.

\bibitem{FM} F. Dillen, M.I. Munteanu,: {\em Constant Angle Surfaces in $\mathbb{H}^2 \times \mathbb{R}$} , Bull. Braz. Math. Soc. Vol. 40(1) ,2009, pp. 85-97. 

\bibitem{GH} G. Ruiz-Hernandez,: {\em Helix, shadow boundary and minimal submanifolds}, Illinois J. Math. Vol. 52(4), 2008, pp. 1385-1397.

\bibitem{PA} P. Cermelly, A.J. Di Scala,:{\em Constant-angle surfaces in liquid cristals}, Philosophical Magazine vol. 87, 2007, pp. 1871-1888.

\bibitem{YTFX} Y. Yin,T. Zhang,F. Yang,X. Qiu, {\em Geometric conditions for fractal supper carbon nonotubes with strict self-similarities} Chaos, Solitons and Fractals, vol. 37(5), 2008, pp. 1257-1266. 


\end{thebibliography}
\end{document}